\documentclass[a4paper]{amsart}
\usepackage{graphicx}
\usepackage{algorithm}
\usepackage{algpseudocode}
\usepackage{amsmath}
\usepackage{amsthm}
\usepackage{amsfonts}
\usepackage{amssymb}
\usepackage{bbm}
\usepackage{bm}
\usepackage{mathrsfs}
\usepackage{mathtools}
\usepackage{thmtools}

\usepackage{tabulary}
\usepackage{booktabs}

\usepackage{fancyref}
\usepackage{hyperref}
\usepackage{listings}
\usepackage{nameref}
\usepackage[section]{placeins}
\usepackage{soul}
\usepackage[normalem]{ulem}
\hypersetup{colorlinks=true, linkcolor=black}
\usepackage{xcolor}

\newtheorem{thm}{Theorem}[section]
\newtheorem{lem}[thm]{Lemma}
\newtheorem{cor}[thm]{Corollary}
\newtheorem{defi}[thm]{Definition}

\newtheorem{prop}[thm]{Proposition}

\declaretheoremstyle[headfont=\normalfont]{normalhead}

\newcommand{\rd}{\mathrm{d}}

\newcommand{\Z}{\mathbb{Z}}

\newcommand{\R}{\mathbb{R}}

\newcommand{\T}{\mathbb{T}}

\newcommand{\Pb}{\mathbb{P}}
\newcommand{\defeq}{\vcentcolon=}
\newcommand{\eqdef}{=\vcentcolon}

\DeclareMathOperator{\dive}{div}

\def\XXint#1#2#3{{\setbox0=\hbox{$#1{#2#3}{\int}$ }
\vcenter{\hbox{$#2#3$ }}\kern-.6\wd0}}

\newcommand{\eps}{\varepsilon}

\graphicspath{ {./HW1/} }
\numberwithin{equation}{section}

\title{A Deterministic Particle Method for the Porous Media Equation}
\author{Amina Amassad and Datong Zhou}

\begin{document}

\maketitle{}

\begin{abstract}
This paper deals with the deterministic particle method for the equation of porous media (with $p=2$). We establish a convergence rate in the Wasserstein-2 distance between the approximate solution of the associated nonlinear transport equation and the solution of the original one.
This seems to be the first quantitative rate for diffusion-velocity particle methods solving diffusive equations and is achieved using a novel commutator estimate for the Wasserstein transport map.

\end{abstract}

\tableofcontents

\section{Introduction}

\subsection{A particle method for the porous media equation}

The purpose of this paper is to shed some light on the convergence properties of a particle method for the Cauchy problem of porous media equation given by
\begin{equation} \label{eq:eq}
\left\{
\begin{aligned}
\partial_t u - \frac{1}{2} \Delta (u^2) & = 0, && \text{in } (0,T) \times \R^d,
\\
u(0,\cdot) & = u_0, && \text{in } \R^d,
\end{aligned} \right.
\end{equation}
where $u \geq 0$ is a density function and $d = 1,2,3$ is the spatial dimension.

The method under consideration is motivated by viewing \eqref{eq:eq} as a transport problem and rewriting
\begin{equation*}
\begin{aligned}
&\partial_t u - \frac{1}{2} \Delta (u^2) = \partial_t u + \dive(V u), &&\text{with } V = - \nabla u.
\end{aligned}
\end{equation*}
To obtain a velocity field with sufficient regularity, the idea is to smooth $u$, and consequently $V$, through convolution with an appropriate kernel, which leads to the following aggregation equation
\begin{equation} \label{eq:num}
\left\{
\begin{aligned}
\partial_t u + \dive(V_\eps u) & = 0, && \text{with } V_\eps = - \nabla R_\eps \star u.
\\
u(0,\cdot) & = u_0, &&
\end{aligned} \right.
\end{equation}
where $R_\eps$ is a mollification kernel of scaling $\eps > 0$. Intuitively, if in some proper sense $R_\eps \to \delta$ as $\eps \to 0$, with $\delta$ being the Dirac mass, then the solution of \eqref{eq:num} converges to a solution of \eqref{eq:eq}.

Let us now demonstrate how \eqref{eq:num} can form the basis of a particle method. For this purpose, let $\{x_i(t)\}_{i=1}^N$ be a  set of $N$ particles, each located at a position $x_i(t) \in \R^d$ for all $t \in (0,T)$. The particle method is realized by evolving the ODE system
\begin{equation} \label{eq:realize}
\begin{aligned}
& \frac{\rd x_i(t)}{\rd t} = - \frac{1}{N} \sum_{j=1}^N \nabla R_\eps\big(x_j(t) - x_i(t) \big), && i = 1,\dots,N.
\end{aligned}
\end{equation}
If we consider the distribution 
\begin{equation*}
\begin{aligned}
u^N(t,x) = \frac{1}{N} \sum_{i = 1}^N \delta\big(x - x_i(t)\big),
\end{aligned}
\end{equation*}
then it can easily be seen that $u^N$ is a solution of \eqref{eq:num} in the sense of distributions on $(0,T) \times \R^d$.
\\

{
In the mathematical literature, particle methods like \eqref{eq:num}, which solve a transport problem with velocity field defined by some heat operator, first appeared in the article \cite{DeMu:90} and have been applied to a variety of problems, for example in \cite{La:99,LaMa:99,LiMa:01}. In particular, \cite{LiMa:01} address the convergence of \eqref{eq:num} to \eqref{eq:eq} as $\eps \to 0$ for periodic solutions. For further review, we refer to \cite{Ch:17} and the references therein. The convergence property of most of these particle methods remain a subject of active research today. Recently, there has been a resurgence of interest in using gradient flow structures. This approach has been remarkably applied to generic porous media equation with exponent $m \neq 2$, as illustrated is studies such as \cite{CaCrPa:19, CaEsSkWu:24, CaEsWu:24}.
}

In this paper, we establish a \emph{rate} (in $\eps$) at which this convergence of \eqref{eq:num} to \eqref{eq:eq} happens.
From the theory developed in \cite{LiMa:01}, it is clear that as the particle number $N$ goes to infinity, $u^N$ the distribution associated with the ODE system \eqref{eq:realize} converges to a solution of \eqref{eq:num} for fixed $\eps > 0$.
A natural question is whether the particle method can give a better approximation of the original equation \eqref{eq:eq} by adapting $\eps \sim N^{-\alpha}$ for an appropriate $\alpha$.
In this paper, we will not be able to have a rigorous discussion on that and confine our attention to the convergence of \eqref{eq:num} in terms of $\eps$.
Moreover, it seems evident to us that the convergence of \eqref{eq:realize} in terms of $\eps$ cannot exceed the convergence of \eqref{eq:num}.

It is worth mentioning that the upcoming analysis also bears relevance for discretizations of \eqref{eq:num} other than \eqref{eq:realize}. For instance, since $V_\eps$ is given by a convolution, this can be computed in $N \log N$ operations through fast Fourier transform. Combined with a standard hyperbolic method for the transport equation, they can provide an explicit numerical method with linear {\sc cfl} condition using $N \log N$ operations at each time step.

\subsection{Gradient flow structure in Wasserstein space}
We address the problem via the optimal transport theory, specifically through the gradient flow structure of the porous media equation in the Wasserstein space. The foundational concepts behind this approach originate from the seminal papers \cite{BeBr:00,Ot:01}. Here, we briefly outline the ideas and postpone the details of optimal transport theory to Section~\ref{subsec:OT}.

The paper considers the periodic solutions in align with \cite{LiMa:01}. To this end, we define the domain of all previously mentioned equations as the torus $\T^d$ and consider the space of probability measure $\Pb(\T^d)$. Everything can be applied to non-negative measures with finite mass, as long as all the measures under consideration have the \emph{same} mass (just divide by the mass).

The weak topology over the space $\Pb(\T^d)$ can be metrized by the $2$-Monge-Kantorovich-Wasserstein distance $W_2$.
Although more straightforward definition of $W_2$ exists, to use its length space and pseudo-Riemannian structure, we prefer the following Benamou-Brenier definition, which is an action-minimizing problem for velocity field $\bm{v}_s$, $s \in [0,1]$ convecting probability measure $\mu_s$, $s \in [0,1]$ from $\mu_0$ to $\mu_1$:
\begin{equation*}
\begin{aligned}
W_2^2(\mu_0, \mu_1) := \inf\left\{\int_0^1 \int |\bm{v}_s(x)|^2 \mu_s(x) \rd x \rd s : \partial_s \mu_s + \dive (\mu_s \bm{v}_s) = 0 \ \text{in distribution} \right\}.
\end{aligned}
\end{equation*}
The action-minimizing curve $\mu_s$, $s \in [0,1]$ serves as the geodesic connecting $\mu_0$, $\mu_1$ and $\bm{v}_s$, $s \in [0,1]$ serves as its tangent vector.
Under proper regularity assumptions, the following formal arguments can be made rigorous: Assume that $\mu_0(t,\cdot)$ and $\mu_1(t,\cdot)$ are probability measures convected by $\bm{w}_{0}(t,\cdot)$, $\bm{w}_{1}(t,\cdot)$, namely
\begin{equation*}
\begin{aligned}
& \partial_t \mu_i + \dive(\bm{w}_{i}\mu_i) = 0, && i = 0,1.
\end{aligned}
\end{equation*}
Then there exists some trajectory $\mu_s(t,\cdot)$ between $\mu_0(t,\cdot)$ and $\mu_1(t,\cdot)$ and a family of velocity fields $\bm{v}_s(t,\cdot)$ such that the $W_2$ distance evolves as
\begin{equation} \label{eq:Wgrad}
\begin{aligned}
\frac{\rd}{\rd t} \bigg[ \frac{1}{2} W_2^2(\mu_0(t,\cdot), \mu_1(t,\cdot)) \bigg] = \int (\bm{v}_1\cdot \bm{w}_{1}\mu_1  - \bm{v}_0\cdot \bm{w}_{0}\mu_0) \ \rd x.
\end{aligned}
\end{equation}
It turns out that the pseudo-Riemannian structure is very well adapted to equations like \eqref{eq:eq}. In particular \eqref{eq:eq} can be interpreted as the gradient flow of functional $E(\mu) = \int \frac{\mu^2}{2} \ \rd x$ in Wasserstein space, which we denote by $\nabla V_{E(\mu)} = \nabla \mu$.
The $W_2$ distance between two gradient flow solutions evolve in a specific manner, namely by substituting $\bm{w}_{0}, \bm{w}_{1}$ in \eqref{eq:Wgrad} by $\nabla V_{E(\mu_0)}$, $\nabla V_{E(\mu_1)}$, one has
\begin{equation} \label{eq:Wgrad_convex}
\begin{aligned}
\int (\bm{v}_1\cdot \nabla V_{E(\mu_1)}\mu_1  - \bm{v}_0\cdot \nabla V_{E(\mu_0)}\mu_0) \ \rd x \leq - \int_0^1\frac{\rd^2}{\rd s^2}E(\mu_s) \ \rd s.
\end{aligned}
\end{equation}
It is well-known today that the choice of $E(\mu)$ above is displacement convex, meaning that $\frac{\rd^2}{\rd s^2}E(\mu_s) \geq 0$ for any action minimizing curve. Consequently, if $\mu_0$ and $\mu_1$ are both solutions of \eqref{eq:eq}, then $\frac{\rd}{\rd t} W_2^2(\mu_0, \mu_1) \leq 0$.

Of course, what we study is the convergence property of \emph{approximate} solutions through \eqref{eq:num}; hence, the perfect non-increasing of $W_2$ distance in time cannot be true.
Nevertheless, we can decompose \eqref{eq:Wgrad} into a main term obeying the gradient flow structure of $E(\mu)$ and some commutator terms quantifying how the evolution of $u^\eps$, the weak solution of \eqref{eq:num} with mollifying kernel $R_\eps$, deviates from \eqref{eq:eq}. (More precisely, we consider the mollification of $u^\eps$ in some intermediate scale.)
The convergence property of $u^\eps \to u$ is then obtained by a commutator estimate, which is the main technical contribution of the paper.

\subsection{Choice of kernel} \label{subsec:kernel}

As will become evident soon, 
the commutator estimate 
depends heavily on the choice of the kernel in \eqref{eq:num}. In this paper, we consider an ``admissible kernel'' belonging to the following case:
\begin{defi} \label{defi:potential}
We say that $R$ is an admissible kernel whenever
\begin{itemize}
\item[i.] $R$ is non-negative in $W^{2,\infty}(\R^d)$ and has bounded second moments. 

\item[ii.] $\int_{\R^d} R(x) \ \rd x = 1$.

\item[iii.] The Fourier transform $\widehat{R}$ is positive and satisfies 
\begin{equation} \label{eq:Rreq_1}
\left\{
\begin{aligned}
\frac{1}{\alpha} \leq & \widehat{R}(\xi) \leq 1, && \forall |\xi| \leq 1,
\\
\frac{a}{|\xi|^{2k}} \leq & \widehat{R}(\xi) \leq \frac{b}{|\xi|^{k}}, && \forall |\xi| > 1,
\end{aligned} \right.
\end{equation}
for some $k > 0$.

\item[iv.] The decomposition $R = R^{\frac{1}{2}} \star R^{\frac{1}{2}}$, defined by $\widehat {R^{\frac{1}{2}}} \defeq (\widehat R)^{\frac{1}{2}}$, satisfies
\begin{equation} \label{eq:Rreq_2}
\begin{aligned}
&|\nabla R^{\frac{1}{2}}(x)| \leq (R^{\frac{1}{2}} \star h)(x), &&\text{a.e. } \ x \in \R^d
\end{aligned}
\end{equation}
for some $h \in \mathcal{M}(\R^d)$, i.e. Borel measure with bounded total variation.

\item[v.] For any $0 < \eta < \eps$, define $R_\eps = \eps^{-d} R(\cdot/\eps)$, $R_\eta = \eta^{-d} R(\cdot/\eta)$ and consider the intermediate scale $L_{\eps,\eta}$ defined as
\begin{equation} \label{eq:Rreq_3}
\begin{aligned}
R_\eps = R_\eta \star L_{\eps,\eta}, \qquad \widehat L_{\eps,\eta} \defeq \frac{\widehat R_\eps}{\widehat R_\eta}.
\end{aligned}
\end{equation}
The kernel $L_{\eps,\eta}$ so defined belongs to $L^1(\R^d)$ and satisfies for some constant $C > 0$,
\begin{equation} \label{eq:Rreq_4}
\begin{aligned}
\int_{\R^d} |x| |L_{\eps,\eta}(x)| \ \rd x \leq C \eps.
\end{aligned}
\end{equation}

\end{itemize}
\end{defi}
Note that due to property iii., $R$ cannot be infinitely smooth and is at best in some appropriate Sobolev spaces.
Moreover, $k$ should be large enough to satisfy the other properties. For example, to expect the embedding
\begin{equation*}
\begin{aligned}
R \in W^{s,2}(\R^d) \subset W^{2,\infty}(\R^d), \quad s > \frac{d}{2}+2,
\end{aligned}
\end{equation*}
a necessary condition is $2k>d+2$.
\\

Extra attention must be given to the 
extension to torus $\T^d$, as it leads to possible ambiguity to define $R_\eps = \eps^{-d}R(\cdot/\eps)$.
To clarify, we start with an admissible kernel $R$ on $\R^d$ to define the corresponding kernel on $\T^d$ through the projection
\begin{equation} \label{eq:R_torus}
\begin{aligned}
& R_\eps^{\T}(x) \defeq \sum_{m \in \Z^d} R_\eps(x + m), && \forall x \in \T^d = [0,1)^d.
\end{aligned}
\end{equation}
and specify our definition of Fourier transform as
\begin{equation*}
\begin{aligned}
\widehat{f}(\xi) \defeq \int_{\R^d} f(x) e^{-2\pi i x \cdot \xi} \ \rd x.
\end{aligned}
\end{equation*}
Then the Fourier series on torus $\T^d$ of such kernel $R_\eps^{\T}$ satisfies
\begin{equation*}
\begin{aligned}
& \widehat{R_\eps^{\T}}(m) = \widehat{R_\eps}(m) = \widehat{R}(\eps m), && \forall \eps > 0, m \in \Z^d.
\end{aligned}
\end{equation*}
We define other kernels, such as $R_\eta^{\T}$ and $L_{\eps,\eta}^{\T}$, through the same projection. All properties in Definition~\ref{defi:potential} are extended in the obvious sense by a careful but straightforward calculation. 
Hence, to simplify the notation, we use the same symbol $R_\eps$ to denote the projection on the torus whenever it does not cause ambiguity.

\subsection{Main results}
The following theorem states our main convergence result.
\begin{thm} \label{thm:main}
Assume initial data $u_0 \in L^\infty(\T^d)$ and an admissible kernel $R$ in the sense of Definition~\ref{defi:potential}. Let $u^\eps \in L^\infty([0,T], L^\infty(\T^d))$ be the weak solution of \eqref{eq:num} with mollifying kernel $R_\eps$ given by \eqref{eq:R_torus}.
Let $u \in L^\infty([0,T], L^\infty(\T^d))$ be the corresponding weak solution of the porous media equation \eqref{eq:eq}. Then there exists a constant $C > 0$, independent of $\eps$, such that
\begin{equation} \label{eq:main_bound}
\begin{aligned}
& \sup_{t \in [0,T]} W_2(u,u^\eps) \leq C \eps^r, && r < \frac{1}{d(4k + 2)},
\end{aligned}
\end{equation}
where $k$ comes from Definition~\ref{defi:potential}. The explicit formula of $C$ depends on $T > 0$, the bounds of $R$ in Definition~\ref{defi:potential} and $\|u_0\|_{L^\infty}$.

\end{thm}
The convergence rate provided by Theorem~\ref{thm:main} declines with increasing regularity of the kernel $R$.
\\

In practice, a rate of convergence in Wasserstein distance may not be the most straightforward. The following corollary, for instance, establishes a convergence result in the $L^2$ sense.
\begin{cor} \label{cor:main}
Under the conditions of Theorem~\ref{thm:main},
\begin{equation*}
\begin{aligned}
& \big\|u - R_\eps^{\frac{1}{2}} \star u^\eps\big\|_{L^2([0,T] \times \T^d)} \leq C \eps^{\frac{2}{2 + d} r},
\end{aligned}
\end{equation*}
where $r$ is given in \eqref{eq:main_bound}.
\end{cor}

In the very special case of dimension $1$, it is possible to get much better rates.
\begin{thm} \label{thm:1d}
Assume that $u_0 \in W^{1,\infty}(\T)$ and that $R$ is a non-negative kernel with $\int_{\R} R(x)\ \rd x = 1$ and $\frac{\rd ^2}{\rd x^2}R(x) \geq 0$ for all $x \neq 0$. Let $u^\eps \in L^\infty([0,T], L^\infty(\T))$ be the weak solution of \eqref{eq:num} with mollifying kernel $R_\eps$ given by \eqref{eq:R_torus}. Then there exists a constant $C > 0$, independent of $\eps$, such that 
\begin{equation*}
\begin{aligned}
& \sup_{t \in [0,T]} W_2(u,u^\eps) \leq C \eps^{\frac{1}{2}}.
\end{aligned}
\end{equation*}
The explicit formula of $C$ depends on $T > 0$, the bounds of $R$ in Definition~\ref{defi:potential} and $\|u_0\|_{W^{1,\infty}}$.

\end{thm}

In Theorem~\ref{thm:1d} the assumptions of Definition~\ref{defi:potential} are completely replaced by a convexity structure. 
Note that since $R$ has to decrease to $0$ at infinity, it is not possible to have $\frac{\rd ^2}{\rd x^2}R(x) \geq 0$ at every point. Instead, the theorem allows $\frac{\rd ^2}{\rd x^2}R(x)$ to contain a minus Dirac mass at $x = 0$. It obviously implies that $R$ cannot be smooth. Still, the structure cannot be generalized in higher dimensions: There does not exist function $R$ s.t. $\nabla^2 R_\eps (x) \geq 0$ (i.e. semi-positive definite) at any $x \neq 0$.
\\

The rest of the paper is structured as follows: Some classical results of \eqref{eq:eq}, \eqref{eq:num} and the preliminaries of optimal transport theory are collected in Section~\ref{subsec:well}~and~\ref{subsec:OT}.
The proof of Theorem~\ref{thm:1d} is given in Section~\ref{subsec:1d}, illustrating the basics of our approach in 1D case.
Section~\ref{sec:proof} is devoted to the more technical proof of Theorem~\ref{thm:main} and Corollary~\ref{cor:main}.

\section{Sketch of the proofs}

\subsection{Well-posedness of equations} \label{subsec:well}
There exists a relatively comprehensive theory concerning the existence, uniqueness and regularity of the porous media equation's solutions.
Obtaining the following result through classical approaches is straightforward. For details, we refer to \cite{BeCrPi:84} and \cite{Va:07}, which address the existence and uniqueness theory for the more general porous media equation $\partial_t u = \frac{1}{m} \Delta ( u^m )$.
\begin{prop} \label{prop:well_posed}
Assume $u_0 \in L^\infty(\T^d)$. Then for any $T > 0$, there exists a unique weak solution $u \in L^\infty([0,T],\ L^\infty(\T^d))$ to \eqref{eq:eq}.
\end{prop}

For a given $R_\eps$, \eqref{eq:num} is a simple aggregation equation with a Lipschitz interaction kernel. Hence, it is immediate to prove that 
\begin{prop} \label{prop:well_posed_num}
Assume $u_0 \in L^\infty(\T^d)$ and $R_\eps \in W^{1,\infty}(\T^d)$. The for any $T > 0$, there exists a unique solution $u^\eps \in L^\infty([0,T],\ L^\infty(\T^d))$ to \eqref{eq:num}. 
\end{prop}
Of course, such $L^\infty$ estimate of \eqref{eq:num} strongly depends on the choice of $R_\eps$. To prove our theorems, we will use quantitative bounds that are independent of $R_\eps$, which will be established in the corresponding sections. 
The role of the $L^\infty$ bounds in our proof is to ensure that the solutions are absolutely continuous with respect to the Lebesgue measure, thereby considerably simplifying the argument in certain instances.

However, in the case of dimension $1$, the following better results are known and will be involved in the quantitative bound in Theorem~\ref{thm:1d}. This further distinguishes dimension $1$ from the general case.
\begin{prop} \label{prop:Lip_1d}
Assume $u_0 \in W^{1,\infty}(\T)$. Then for any $T > 0$, the solution $u$ to \eqref{eq:eq} belongs to $L^\infty([0,T],\ L^\infty(\T^d))$.
\end{prop}

\subsection{Tools from optimal transport} \label{subsec:OT}
The aim of this subsection is to recall some preliminary results regarding optimal transport theory that are available, for instance, in \cite{BeBr:00, Ot:01}. We refer to \cite{AmGiSa:05, Vi:09} for detailed discussions.

First, let us revisit several equivalent definitions of $W_2$ distance on $\T^d$. Let $|x-y|_{\T^{d}}$ denote the projection of the Euclidian norm of $\R^d$ on $\T^d$.
The Kantorovich formulation of $W_2$ is given by
\begin{equation*}
\begin{aligned}
W_2^2(\mu_0, \mu_1) = \inf\bigg\{ & \int_0^1 |x-y|_{\T^{d}}^2 \ \pi(\rd x, \rd y) :
\\
& \pi \in \Pb(\T^d \times \T^d), \int_{\T^d} \pi(x,\rd y) = \mu_0(x), \int_{\T^d} \pi(\rd x, y) = \mu_1(y),
\bigg\}.
\end{aligned}
\end{equation*}
This can be reduced to the following Monge formulation when $\mu_0$ is absolutely continuous with respect to the Lebesgue measure:
\begin{equation*}
\begin{aligned}
W_2^2(\mu_0, \mu_1) = \inf\bigg\{ & \int_0^1 |x-T(x)|_{\T^{d}}^2 \ \mu_0(\rd x) : T_\#(\mu_0) = \mu_1 \bigg\}.
\end{aligned}
\end{equation*}

Let us now return to the aforementioned Benamou-Brenier formulation
\begin{equation*}
\begin{aligned}
W_2^2(\mu_0, \mu_1) = \inf\bigg\{\int_0^1 \int |\bm{v}_s(x)|^2 \mu_s(x) \rd x \rd s : \partial_s \mu_s + \dive (\mu_s \bm{v}_s) = 0 \ \text{in distribution} \bigg\}.
\end{aligned}
\end{equation*}
For general probability measure $\mu_0, \mu_1$, the optimal $\mu_s$ and $\bm{v}_s$ can be very singular and should be understood in the sense of convex functional for measures. However, on compact domains such as $\T^d$ when at least one of $\mu_0$ is absolutely continuous w.r.t. Lebesgue measure, the action-minimizing convection $\bm{v}_s$ is much more regular. In particular it can be written as $\bm{v}_s = \nabla \psi_s$, where $\nabla \psi_s$ solving the Hamilton-Jacobi equation
\begin{equation*}
\begin{aligned}
& \partial_s \psi_s + \frac{|\nabla \psi_s|^2}{2} = 0, && s \in [0,1].
\end{aligned}
\end{equation*}
The potential $\psi_s$ satisfies the so-called $d^2/2$-convex condition. On the torus $\T^d$ it implies that by extending $\psi_s$ to a periodic function on $\R^d$, the function $\phi_s = \psi_s + \frac{|x|^2}{2}$ is convex. Then $T(x) = \nabla \phi_0(x) = \nabla \psi_0(x) + x$ give the convex optimal Monge map pushing $\mu_0$ to $\mu_1$.
Indeed, on $\T^d$, the travel distance of an optimal map cannot exceed $\sqrt{d} / 2$. Since $\phi_s$ is convex, this implies some BV regularity of $\bm{v}_s$.
In summary, we have
\begin{lem} \label{lem:BV}
Let $\mu_0, \mu_1 \in \Pb(\T^d)$. Assume that $\mu_0$ or $\mu_1$ is absolutely continuous w.r.t. the Lebesgue measure. Then there exists a constant $C$ such that an action-minimizing $\bm{v}_s, s \in [0,1]$ in the Benamou-Brenier formulation satisfies
\begin{equation*}
\begin{aligned}
\sup_{s} \big( \|\bm{v}_s\|_{L^\infty(\T^d)} + \|\bm{v}_s\|_{BV(\T^d)} \big) \leq C.
\end{aligned}
\end{equation*}
\end{lem}
Such a priori regularity of $\bm{v}_s$ will be a key component in our commutator estimate for the $W_2^2$ evolution, heuristically previewed by \eqref{eq:Wgrad} and rigorous detailed as follows:
\begin{lem} \label{lem:W2diff}
Let $\mu_0, \mu_1 \in L^\infty([0,T], \Pb(\T^d))$, absolutely continuous w.r.t. the Lebesgue measure.
Assume that they solve the convection equations
\begin{equation*}
\begin{aligned}
& \partial_t \mu_i + \dive(\bm{w}_{i}\mu_i) = 0, && i = 0,1,
\end{aligned}
\end{equation*}
in distribution, with
\begin{equation*}
\begin{aligned}
& \int_0^T \int |\bm{w}_i(t,x)|^2 \mu_i(t,x) \rd x \rd t < \infty, && i = 0,1.
\end{aligned}
\end{equation*}
Then there exists an action-minimizing $\mu_s$, $\bm{v}_s$, $s \in [0,1]$ in the Benamou-Brenier formulation, such that for a.e. $t \in [0,T]$
\begin{equation} \label{eq:Wgrad_ineq}
\begin{aligned}
\frac{\rd^+}{\rd t} \bigg[ \frac{1}{2} W_2^2(\mu_0(t,\cdot), \mu_1(t,\cdot)) \bigg] \leq \int (\bm{v}_1\cdot \bm{w}_{1}\mu_1  - \bm{v}_0\cdot \bm{w}_{0}\mu_0) \ \rd x.
\end{aligned}
\end{equation}
Moreover by Lemma~\ref{lem:BV}, there exists some constant $C > 0$ such that
\begin{equation*}
\begin{aligned}
\sup_{s,t} \big( \|\bm{v}_s(t,\cdot)\|_{L^\infty(\T^d)} + \|\bm{v}_s(t,\cdot)\|_{BV(\T^d)} \big) \leq C.
\end{aligned}
\end{equation*}

\end{lem}
The inequality in \eqref{eq:Wgrad_ineq} can be made equality for derivative by an in-depth investigation of the regularity of flows. However, for the purpose of obtaining stability of gradient flow, the one-sided result is all we need. 

Many diffusion equations like \eqref{eq:eq} have a gradient flow structure adapted to the above pseudo-Riemannian structure of $W_2$ distance.
Consider a functional of form
\begin{equation*}
\begin{aligned}
E(\mu) = \int \Phi(\mu(x)) \ \rd x,
\end{aligned}
\end{equation*}
One may formally define the gradient of $E(\mu)$ with respect to the $W_2$ metric in the following framework, called Otto's calculus since \cite{Ot:01}: For any trajectory $\mu_s$ such that
\begin{equation} \label{eq:Wcurve}
\begin{aligned}
\partial_t \mu_t + \dive (\mu_t \bm{w}_t) = 0,
\end{aligned}
\end{equation}
the time derivative of $E(\mu)$ should be the inner product of gradient of $E(\mu)$ and the field $\bm{w}_t$ (as tangent vector of $\mu_t$). By the identity
\begin{equation*}
\begin{aligned}
\frac{\rd}{\rd t} E(\mu_t) \Big|_{t = 0} = - \int \Phi'(\mu_0) \dive (\mu_0 \bm{w}_0) \ \rd x = \int \nabla( \Phi'(\mu_0)) (\mu_0 \bm{w}_0) \ \rd x,
\end{aligned}
\end{equation*}
it is easy to see that the formal gradient vector at $\mu_0$ is $\nabla( \Phi'(\mu_0))$. A function then solves the gradient flow with respect to $E$ iff it solves \eqref{eq:Wcurve} with
\begin{equation*}
\begin{aligned}
\bm{w}_t = \nabla V_{E(\mu)} = - \nabla( \Phi'(\mu)).
\end{aligned}
\end{equation*}
For instance, if $\Phi(\mu)=\frac{1}{2}\mu^2$, then the gradient flow is exactly the porous media equation \eqref{eq:eq}. The most classical case is $\Phi(\mu)=\mu\,\log \mu$, where the corresponding gradient flow is simply the heat equation.
One way to rigorously justify the above formal arguments involves the subdifferential on metric spaces; see \cite{AmGiSa:05}.

The above formal Otto's calculus also suggests the $W_2$ distance of two measures convected by gradient flows evolves in a specific way: Let $\mu_s$, $\bm{v}_s$, $s \in [0,1]$ be an action-minimizing path in the sense of Lemma~\ref{lem:W2diff},
\begin{equation*} \label{eq:Wgrad_flow}
\begin{aligned}
& \int (\bm{v}_1\cdot \nabla V_{E(\mu_1)} \mu_1  - \bm{v}_0\cdot \nabla V_{E(\mu_0)} \mu_0) \ \rd x
\\
& \quad = - \bigg( \frac{\rd}{\rd s} E(\mu_s) \Big|_{s = 1} \bigg) + \bigg( \frac{\rd}{\rd s} E(\mu_s) \Big|_{s = 0} \bigg)
\end{aligned}
\end{equation*}
In particular, if for \emph{any} action-minimizing $\mu_s$, the function $s \mapsto E(\mu_s)$ is convex, one says that $E$ is \emph{displacement convex} in $W_2$ distance.
As a consequence, such gradient flow is contracting for the $W_2$ distance:
\begin{lem} \label{lem:Wgrad_contract}
Let $\mu_0, \mu_1 \in \Pb(\T^d)$ and that $\nabla V_{E(\mu_0)}$, $\nabla V_{E(\mu_1)}$ are the gradient flow of a function $E$ being displacement convex in $W_2$ distance. Denote $\mu_s$, $\bm{v}_s$, $s \in [0,1]$ the action-minimizing path between $\mu_0, \mu_1$ in the sense of Lemma~\ref{lem:W2diff}. Then
\begin{equation*} \label{eq:Wgrad_contract}
\begin{aligned}
\int (\bm{v}_1\cdot \nabla V_{E(\mu_1)} \mu_1  - \bm{v}_0\cdot \nabla V_{E(\mu_0)} \mu_0) \ \rd x \leq 0.
\end{aligned}
\end{equation*}

\end{lem}
The most classical example of displacement convexity is $E(\mu) = \int \Phi(\mu(x)) \rd x$ with $\Phi: \R_+ \to \R$ being a convex function.
Hence, the solutions of equation \eqref{eq:eq} are contracting. 

In addition, the gradient flow structure applies to the aggregation equation \eqref{eq:num}, which was first studied in \cite{CaMcVi:03}.
The corresponding functional is of form
\begin{equation*}
\begin{aligned}
E_\eps(\mu) = \frac{1}{2} \int \mu \big( R_\eps \star \mu \big) \ \rd x.
\end{aligned}
\end{equation*}
However, $E_\eps(\mu)$ is displacement convex only if $\nabla^2 R_\eps (x) \geq 0$ (i.e. semi-positive definite) for any $x \neq 0$. This cannot be true if $\int |R_\eps(x)| \ \rd x < \infty$, except in dimension $1$.
This displacement convexity allows a much more straightforward analysis in dimension $1$.

\subsection{Proof of Theorem~\ref{thm:1d}} \label{subsec:1d}
We conclude the section by the proof for the dimensional $1$ case. 

The key commutator estimate is made through a direction somewhat opposite to the one in the general proof.
Namely, it quantifies how the evolution of $u$ deviates from \eqref{eq:num} instead of how the evolution of $u^\eps$ deviates from \eqref{eq:eq}.
This leads to a better rate in dimension $1$ exclusively since, as explained previously, one cannot expect \eqref{eq:num} to be contracting in higher dimensions.

\begin{proof} [Proof of Theorem~\ref{thm:1d}]
Throughout the proof, the kernel $R$ and scaling factor $\eps$ will be fixed. For simplicity, we denote $\tilde u := u^\eps$, a solution to \eqref{eq:num}.

Apply Lemma~\ref{lem:W2diff} to $u$ and $\tilde u$ to find some $\bm{v}_s(t,x)$ such that
\begin{equation*}
\begin{aligned}
\frac{\rd^+}{\rd t} \bigg[ \frac{1}{2} W_2^2(u,\tilde u) \bigg] \leq \int \Big( \bm{v}_1\cdot ( - \nabla R_\eps \star \tilde u) \tilde u  - \bm{v}_0\cdot (- \nabla u) u \Big) \ \rd x \eqdef D_{u, \tilde u}.
\end{aligned}
\end{equation*}
By adding and subtracting,
\begin{equation*}
\begin{aligned}
D_{u, \tilde u} & = \int \Big( \bm{v}_1\cdot ( - \nabla R_\eps \star \tilde u) \tilde u  - \bm{v}_0\cdot (- \nabla R_\eps \star u) u \Big) \ \rd x  -  \int \bm{v}_0\cdot (\nabla R_\eps \star u - \nabla u) u \ \rd x
\\
& \eqdef G_{u, \tilde u} - C_{u, \tilde u}.
\end{aligned}
\end{equation*}
The functional $E_\eps$ is displacement convex since here $\frac{\rd^2}{\rd x^2} R_\eps \geq 0$ for $x \neq 0$, thus by Lemma~\ref{eq:Wgrad_contract}, 
\begin{equation*}
\begin{aligned}
G_{u, \tilde u} \leq 0.
\end{aligned}
\end{equation*}

By integration by part,
\begin{equation*}
\begin{aligned}
C_{u, \tilde u} = - \int \dive(u \ \bm{v}_0) ( R_\eps \star u - u) \ \rd x
\end{aligned}
\end{equation*}
Now $u$ is bounded in $W^{1,\infty}$ by Proposition~\ref{prop:Lip_1d} and $\bm{v}_0$ is bounded uniformly in $\eps$ in $L^\infty \cap BV$ by Lemma~\ref{lem:W2diff}. Therefore $\dive(u \ \bm{v}_0)$ is a measure of bounded total variance. Consequently
\begin{equation*}
\begin{aligned}
|C_{u, \tilde u}| \leq C \|R_\eps \star u - u\|_{L^\infty} \leq C \eps,
\end{aligned}
\end{equation*}
as again $u$ is Lipschitz by Proposition~\ref{prop:Lip_1d}.
The constant $C$ may vary through the inequalities, and for the righthand side, it depends on the choice of $R$, $\|u_0\|_{W^{1,\infty}}$ and $T > 0$. (The last two essentially control the a priori $W^{1,\infty}$-bound of $u$.)

Integrate over time to get
\begin{equation*}
\begin{aligned}
\frac{1}{2} W_2^2(u,\tilde u) \leq \int_0^T |C_{u, \tilde u}| \ \rd t \leq C \eps.
\end{aligned}
\end{equation*}
We conclude the assertion in the theorem by taking the square root on both sides.
\end{proof}

\section{Proof in general dimensions} \label{sec:proof}
Like in the proof of Theorem~\ref{thm:1d}, from now on, we denote $\tilde u := u_\eps$ a solution to \eqref{eq:num}.
Moreover, we denote its mollifications by $\tilde u_\eps=R_\eps\star \tilde u$ and $\tilde u_\eta=R_\eta \star \tilde u$.

Again, we initiate the proof by applying Lemma~\ref{lem:W2diff} to $u$ and $\tilde u$ to find some $\bm{v}_s(t,x)$ such that
\begin{equation*}
\begin{aligned}
\frac{\rd^+}{\rd t} \bigg[ \frac{1}{2} W_2^2(u,\tilde u) \bigg] \leq \int \Big( \bm{v}_1\cdot ( - \nabla R_\eps \star \tilde u) \tilde u  - \bm{v}_0\cdot (- \nabla u) u \Big) \ \rd x \eqdef D_{u, \tilde u}.
\end{aligned}
\end{equation*}
However, the commutator estimate is much more complicated, and we will take several steps to make it.

\subsection{Fundamental a priori estimates}
The following stability estimate follows in a straightforward manner using the decomposition $R_\eps = R_\eps^{\frac{1}{2}} \star R_\eps^{\frac{1}{2}}$.
\begin{lem} \label{lem:energy}
Let $R$ be an admissible kernel in the sense of Definition~\ref{defi:potential} and set $R_\eps(\cdot) = \eps^{-d} R(\cdot / \eps)$. Given initial data $\tilde u_0$ such that $R_\eps^{\frac{1}{2}} \star \tilde u_0 \in L^2$ and $\tilde u_0 \log \tilde u_0 \in L^1$, there exists a solution $\tilde u$ of the aggregation equation \eqref{eq:num} satisfying the following energy estimates
\begin{equation*}
\begin{aligned}
\sup_{t} \bigg( \Big\|R_\eps^{\frac{1}{2}} \star \tilde u\Big\|_{L^2}^2 + \int_0^t \int \tilde u |\nabla \tilde u_\eps|^2 \ \rd x \rd \tau \bigg) & \leq \Big\|R_\eps^{\frac{1}{2}} \star \tilde u_0\Big\|_{L^2}^2,
\\
\sup_{t} \bigg( \int \tilde u \log \tilde u \ \rd x + \int_0^t \int \big| \nabla R_\eps^{\frac{1}{2}} \star \tilde u \big|^2 \ \rd x \rd \tau \bigg) & \leq \int \tilde u_0 \log \tilde u_0 \ \rd x.
\end{aligned}
\end{equation*}

\end{lem}

\begin{proof}
Using \eqref{eq:num} and the symmetry of $R_\eps$, we have that 
\begin{equation*}
\begin{aligned}
\frac{\rd}{\rd t} \frac{1}{2}\Big\|R_\eps^{\frac{1}{2}} \star \tilde u\Big\|_{L^2}^2
= \int \partial_t \tilde u \  (R_\eps \star \tilde u)
= \int -\dive \big( \tilde u (-\nabla R_\eps \star \tilde u) \big) \  (R_\eps \star \tilde u) \ \rd x.
\end{aligned}
\end{equation*}
Integrate by parts and recall that $\tilde u_\eps = R_\eps \star \tilde u$, we conclude that 
\begin{equation*}
\begin{aligned}
\frac{\rd}{\rd t} \frac{1}{2}\Big\|R_\eps^{\frac{1}{2}} \star \tilde u\Big\|_{L^2}^2 = - \int  \tilde u |\nabla \tilde u_\eps|^2 \ \rd x,
\end{aligned}
\end{equation*}
from which the first estimate in the lemma follows.

Next, we calculate that
\begin{equation*}
\begin{aligned}
\frac{\rd}{\rd t} \int \tilde u \log \tilde u \ \rd x = \int \partial_t \tilde u \ (1 + \log \tilde u) \ \rd x = \int -\dive \big( \tilde u (-\nabla R_\eps \star \tilde u) \big) \  (1 + \log \tilde u) \ \rd x.
\end{aligned}
\end{equation*}
Integrate by parts and use the symmetry of $R_\eps$, we have that
\begin{equation*}
\begin{aligned}
\frac{\rd}{\rd t} \int \tilde u \log \tilde u \ \rd x = \int  (-\nabla R_\eps \star \tilde u) \  \nabla \tilde u \ \rd x = - \int \big| \nabla R_\eps^{\frac{1}{2}} \star \tilde u \big|^2 \ \rd x.
\end{aligned}
\end{equation*}
Integrating in time gives the second estimate in the lemma.
\end{proof}
Note that on $\T^d$, the integration $\int \tilde u \log \tilde u \ \rd x$ is uniformly bounded from below and $\int \tilde u_0 \log \tilde u_0 \ \rd x$ is bounded from above by $\log( \|u_0\|_{L^\infty})$. Hence a direct consequence of Lemma~\ref{lem:energy} is the a priori boundedness uniform in $\eps > 0$:
\begin{equation*}
\begin{aligned}
\int_0^t \Big\|R_\eps^{\frac{1}{2}} \star \tilde u\Big\|_{H^1}^2 \ \rd \tau < C.
\end{aligned}
\end{equation*}
However, it does not give us enough Sobolev regularity of $\tilde u$ without mollification. Therefore, estimating a naive commutator $C_{u, \tilde u}$ of form
\begin{equation*}
\begin{aligned}
D_{u, \tilde u}
& = \int \Big( \bm{v}_1\cdot ( - \nabla \tilde u) \tilde u  - \bm{v}_0\cdot (- \nabla u) u \Big) \ \rd x - \int \Big( \bm{v}_1\cdot ( \nabla R_\eps \star \tilde u - \nabla \tilde u) \tilde u \Big) \ \rd x
\\
& \eqdef G_{u, \tilde u} - C_{u, \tilde u}
\end{aligned}
\end{equation*}
is still out of our reach. Instead, we introduce an intermediate scale $R_\eta$ that will be specified later. The mollified $\tilde u_\eta$ inherits some regularity from the above estimate; hence, it is a good starting point for our commutator estimate.

\subsection{Kernel decomposition}

In this subsection, we establish the estimates involving intermediate scale mollification $R_\eta$, where $\eta < \eps$.
\begin{lem} \label{lem:intermediate_1}
Let $R$ be an admissible kernel in the sense of Definition~\ref{defi:potential} and $0 < \eta < \eps$. Then for all $f \in H^{-1}$, 
\begin{equation*}
\begin{aligned}
\|R_\eta \star f\|_{L^2} \leq C \bigg( \frac{\eps}{\eta} \bigg)^k \Big\|R_\eps^{\frac{1}{2}} \star f\Big\|_{L^2},
\end{aligned}
\end{equation*}
where $C = \big(\alpha + \frac{1+b^2}{a})^{\frac{1}{2}}$ (see Definition~\ref{defi:potential}).

\end{lem}

\begin{proof}
We begin with the proof for $f \in L^2$, so there is no issue when applying the Fourier transform.
The following estimate actually holds for both $f \in L^2(\R^d)$ and $f \in L^2(\T^d)$, provided that the integration is understood as a summation over lattice $\xi \in \Z^d$ in the case of $\T^d$.
\begin{equation*}
\begin{aligned}
& \|R_\eta \star f\|_{L^2}^2
\\
& \quad = \int \big|\widehat R_\eta \widehat f \big|^2 \mathbbm{1}_{\{0 \leq |\xi| < \frac{1}{\eps}\}} \ \rd \xi + \int \big|\widehat R_\eta \widehat f \big|^2 \mathbbm{1}_{\{\frac{1}{\eps} \leq |\xi| < \frac{1}{\eta}\}} \ \rd \xi + \int \big|\widehat R_\eta \widehat f \big|^2 \mathbbm{1}_{\{\frac{1}{\eta} \leq |\xi|\}} \ \rd \xi
\\
& \quad \eqdef F_1 + F_2 + F_3.
\end{aligned}
\end{equation*}
Recall the property \textnormal{iii.} in Definition~\ref{defi:potential}
\begin{subequations}
\begin{align}
\frac{1}{\alpha} \leq & \widehat{R}(\xi) \leq 1, && \forall |\xi| \leq 1, \label{eq:iiia}
\\
\quad \quad \quad \quad \quad \quad \quad \quad \frac{a}{|\xi|^{2k}} \leq & \widehat{R}(\xi) \leq \frac{b}{|\xi|^{k}}, && \forall |\xi| > 1, \label{eq:iiib}
\end{align}
\end{subequations}
we have the follow estimates: Using \eqref{eq:iiia} it is straightforward that
\begin{equation*}
\begin{aligned}
F_1 = \int_{0 \leq |\xi| < \frac{1}{\eps}} \big|\widehat R_\eta \widehat f \big|^2 \ \rd \xi
=
\int_{0 \leq |\xi| < \frac{1}{\eps}} \frac{\widehat R_\eta^2}{\widehat R_\eps} \big|\widehat R_\eps^{\frac{1}{2}} \widehat f \big|^2 \ \rd \xi
\leq
\int_{0 \leq |\xi| < \frac{1}{\eps}} \alpha \big|\widehat R_\eps^{\frac{1}{2}} \widehat f \big|^2 \ \rd \xi.
\end{aligned}
\end{equation*}
By using the second inequality in \eqref{eq:iiia} then the first inequality in \eqref{eq:iiib}, we calculate
\begin{equation*}
\begin{aligned}
F_2 = \int_{\frac{1}{\eps} \leq |\xi| < \frac{1}{\eta}} \big|\widehat R_\eta \widehat f \big|^2 \ \rd \xi & \leq \int_{\frac{1}{\eps}
\leq
|\xi| < \frac{1}{\eta}} \big| \widehat f \big|^2 \ \rd \xi
\\
& = \int_{\frac{1}{\eps} \leq |\xi| < \frac{1}{\eta}} \frac{1}{\widehat R_\eps} \big|\widehat R_\eps^{\frac{1}{2}} \widehat f \big|^2 \ \rd \xi
\leq
\int_{\frac{1}{\eps} \leq |\xi| < \frac{1}{\eta}} \frac{1}{a} \bigg( \frac{\eps}{\eta} \bigg)^{2k} \big|\widehat R_\eps^{\frac{1}{2}} \widehat f \big|^2 \ \rd \xi.
\end{aligned}
\end{equation*}
Last, using \eqref{eq:iiib}, we observe
\begin{equation*}
\begin{aligned}
F_3 = \int_{\frac{1}{\eta} \leq |\xi|} \big|\widehat R_\eta \widehat f \big|^2 \ \rd \xi
=
\int_{\frac{1}{\eta} \leq |\xi|} \frac{\widehat R_\eta^2}{\widehat R_\eps} \big|\widehat R_\eps^{\frac{1}{2}} \widehat f \big|^2 \ \rd \xi
\leq
\int_{\frac{1}{\eta} \leq |\xi|} \frac{b^2}{a} \bigg( \frac{\eps}{\eta} \bigg)^{2k} \big|\widehat R_\eps^{\frac{1}{2}} \widehat f \big|^2 \ \rd \xi.
\end{aligned}
\end{equation*}
Combining the three bounds, we conclude the lemma for $f \in L^2$.

The estimate extends to $f \in H^{-1}$ by a density argument. WLOG assume $f = \partial_{x_1} \varphi$ for some $\varphi \in L^2$. By choosing $\tilde \varphi \in H^1$ approximating $\varphi$, we have
\begin{equation*}
\begin{aligned}
& \|R_\eta \star \partial_{x_1} \varphi \|_{L^2}
\\
& \quad \leq
\|R_\eta \star \partial_{x_1} \tilde \varphi \|_{L^2} + \|\partial_{x_1} R_\eta \star  (\varphi - \tilde \varphi) \|_{L^2}
\\
& \quad \leq C \bigg( \frac{\eps}{\eta} \bigg)^k \Big\|R_\eps^{\frac{1}{2}} \star \partial_{x_1} \tilde \varphi \Big\|_{L^2} + \|\partial_{x_1} R_\eta \star  (\varphi - \tilde \varphi) \|_{L^2}
\\
& \quad \leq C \bigg( \frac{\eps}{\eta} \bigg)^k \Big\|R_\eps^{\frac{1}{2}} \star \partial_{x_1} \varphi \Big\|_{L^2} + C \bigg( \frac{\eps}{\eta} \bigg)^k \Big\|\partial_{x_1} R_\eps^{\frac{1}{2}} \star (\varphi - \tilde \varphi) \Big\|_{L^2} + \|\partial_{x_1} R_\eta \star  (\varphi - \tilde \varphi) \|_{L^2}.
\end{aligned}
\end{equation*}
The remainder terms vanish as $\tilde \varphi \to \varphi \in L^2$.

\end{proof}

The following inverse estimate on derivative relies on the non-negative assumption:
\begin{lem} \label{lem:intermediate_2}
Let $R$ be an admissible kernel in the sense of Definition~\ref{defi:potential} and $0 < \eta < \eps$. Then for all non-negative $f \in L^2$, 
\begin{equation*}
\begin{aligned}
\Big\| \big| \nabla R_\eps^{\frac{1}{2}} \big| \star f\Big\|_{L^2} \leq C \bigg( \frac{1}{\eps} \bigg) \Big\| R_\eps^{\frac{1}{2}} \star f\Big\|_{L^2}
\end{aligned}
\end{equation*}
where $C = \|h\|_{TV}$ (see Definition~\ref{defi:potential}).
\end{lem}
\begin{proof}
Recall property iv. of Definition~\ref{defi:potential}, $|\nabla R^{\frac{1}{2}}(x)| \leq (R^{\frac{1}{2}} \star h)(x)$, a.e. $x \in \R^d$, it is straightforward to verify that for a.e. $x \in \R^d$,
\begin{equation*}
\begin{aligned}
& (R_\eps^{\frac{1}{2}} \star h_\eps)(x)
\\
& \quad = \int_{\R^d} \frac{1}{\eps^d} R^{\frac{1}{2}}\bigg( \frac{z - x}{\eps} \bigg) \frac{1}{\eps^d} h \bigg( \frac{z}{\eps} \bigg) \ \rd z
\\
& \quad = \frac{1}{\eps^d} (R^{\frac{1}{2}} \star h)(x/\eps)
\\
& \quad \geq \frac{1}{\eps^d} \big| \nabla R^{\frac{1}{2}} (x/\eps) \big|.
\end{aligned}
\end{equation*}
By the chain rule, we further have $(R_\eps^{\frac{1}{2}} \star h_\eps)(x) \geq \eps \big| \nabla R_\eps^{\frac{1}{2}} (x) \big|$, $\forall x \in \R^d$. This also extends to $\T^d$ by choosing proper projection as specified in Section~\ref{subsec:kernel}. By non-negativity of $f$,
\begin{equation*}
\begin{aligned}
\Big\| \big| \nabla R_\eps^{\frac{1}{2}} \big| \star f\Big\|_{L^2}
\leq
\frac{1}{\eps} \Big\|R_\eps^{\frac{1}{2}} \star h_\eps \star f \Big\|_{L^2}
\leq
C \bigg( \frac{1}{\eps} \bigg) \Big\| R_\eps^{\frac{1}{2}} \star f\Big\|_{L^2}.
\end{aligned}
\end{equation*}
\end{proof}

The two lemmas are applied in our commutator estimate for $f = \nabla \tilde u$ and $f = \tilde u$, respectively. This is where the $\eps$-dependent $L^\infty$-estimate of $\tilde u$ gets involved to fulfill the non-quantitative part of the regularity requirement.

\subsection{Key commutator estimate}
We now establish the commutator estimate for the intermediate scale $\tilde u_\eta = R_\eta \star \tilde u$.
The following proposition will be our main result of this subsection:
\begin{prop} \label{prop:grad_inter}
Let $u_0 \in L^\infty(\T^d)$, $\tilde u = u^\eps$ solve \eqref{eq:num} and $u$ solve \eqref{eq:eq}. Then for any $\eps > 0$, any $\gamma < 1 + (d(2k+1))^{-1}$, take $\eta = \eps^{\gamma}$, there exists a constant $C > 0$, independent of $\eps$, such that
\begin{equation*}
\begin{aligned}
&\sup_t W_2^2(u, \tilde u_\eta) \leq C \eps^r, && r = \gamma - 1 < \frac{1}{d(2k+1)},
\end{aligned}
\end{equation*}
where $k$ is given by Definition~\ref{defi:potential}.
\end{prop}
\begin{proof}
Notice that 
\begin{equation} \label{eq:num_conv}
\begin{aligned}
\partial_t \tilde u_\eta = \dive \big( R_\eta \star (\tilde u \ \nabla \tilde u_\eps) \big) = \dive \bigg( \frac{R_\eta \star (\tilde u \ \nabla \tilde u_\eps)}{\tilde u_\eta} \tilde u_\eta \bigg).
\end{aligned}
\end{equation}
Apply Lemma~\ref{lem:W2diff} to $u$ and $\tilde u_\eta$ to find some $\bm{v}_s(t,x)$ such that
\begin{equation*}
\begin{aligned}
\frac{\rd^+}{\rd t} \bigg[ \frac{1}{2} W_2^2(u,\tilde u_\eta) \bigg] \leq \int \bigg( \bm{v}_1\cdot \big( - R_\eta \star (\tilde u \ \nabla \tilde u_\eps) \big) - \bm{v}_0\cdot (- \nabla u) u \bigg) \ \rd x \eqdef D_{u, \tilde u_\eta}.
\end{aligned}
\end{equation*}
We construct their commutator terms by the following editing:
\begin{equation*}
\begin{aligned}
& \int \bm{v}_1\cdot \big( R_\eta \star (\tilde u \ \nabla \tilde u_\eps) \big) \ \rd x
\\
& \quad \eqdef \int \bm{v}_1\cdot (R_\eta \star \tilde u) \ \nabla \tilde u_\eps \ \rd x + C_{u, \tilde u_\eta}^{(1)}
\\
& \quad = \int \bm{v}_1\cdot \tilde u_\eta \ (\nabla L_{\eps,\eta} \star R_\eta \star \tilde u ) \ \rd x + C_{u, \tilde u_\eta}^{(1)}
\\
& \quad \eqdef \int \bm{v}_1\cdot ( \tilde u_\eta \ \nabla \tilde u_\eta ) \star L_{\eps,\eta} \ \rd x + C_{u, \tilde u_\eta}^{(1)} + C_{u, \tilde u_\eta}^{(2)}
\\
& \quad \eqdef \int ( \bm{v}_1\cdot \tilde u_\eta \ \nabla \tilde u_\eta ) \star L_{\eps,\eta} \ \rd x + C_{u, \tilde u_\eta}^{(1)} + C_{u, \tilde u_\eta}^{(2)} + C_{u, \tilde u_\eta}^{(3)}.
\end{aligned}
\end{equation*}
The exact formulation of the commutators will be reviewed later.
Observe that the last integral 
\begin{equation*}
\begin{aligned}
\int ( \bm{v}_1\cdot \tilde u_\eta \ \nabla \tilde u_\eta ) \star L_{\eps,\eta} \ \rd x = \int \bm{v}_1\cdot \tilde u_\eta \ \nabla \tilde u_\eta \ \rd x,
\end{aligned}
\end{equation*}
which follows the gradient flow structure of \eqref{eq:eq}. Therefore by Lemma~\ref{eq:Wgrad_contract}
\begin{equation*}
\begin{aligned}
G_{u, \tilde u_\eta} \defeq \int \bigg( \bm{v}_1\cdot (- \nabla \tilde u_\eta) \tilde u_\eta - \bm{v}_0\cdot (- \nabla u) u \bigg) \ \rd x \leq 0,
\end{aligned}
\end{equation*}
and $D_{u, \tilde u_\eta}$ can be rewritten as
\begin{equation*}
\begin{aligned}
D_{u, \tilde u_\eta} = G_{u, \tilde u_\eta} - ( C_{u, \tilde u_\eta}^{(1)} + C_{u, \tilde u_\eta}^{(2)} + C_{u, \tilde u_\eta}^{(3)} ).
\end{aligned}
\end{equation*}
It remains to bound the three commutator terms. \\

\paragraph{\emph{Bound of $C_{u, \tilde u_\eta}^{(1)}$}} By applying the decomposition $R_\eps = R_\eps^{\frac{1}{2}} \star R_\eps^{\frac{1}{2}}$, we deduce
\begin{equation*}
\begin{aligned}
C_{u, \tilde u_\eta}^{(1)} & \defeq
\int \bm{v}_1\cdot \big( R_\eta \star (\tilde u \ \nabla \tilde u_\eps) \big) \ \rd x - \int \bm{v}_1\cdot (R_\eta \star \tilde u) \ \nabla \tilde u_\eps \ \rd x
\\
& = \iint \bm{v}_1(x) \cdot R_\eta(x - y) (\nabla \tilde u_\eps(y) - \nabla \tilde u_\eps(x)) \ \tilde u(y) \ \rd y \rd x
\\
& = \iiint \bm{v}_1(x) \cdot R_\eta(x - y) \Big[ (R_\eps^{\frac{1}{2}} \star \tilde u)(y - z) - (R_\eps^{\frac{1}{2}} \star \tilde u)(x - z)\Big] \ \nabla R_\eps^{\frac{1}{2}}(z) \tilde u(y) \ \rd z \rd y \rd x
\end{aligned}
\end{equation*}
Notice that all terms are non-negative except $\bm{v}_1$ and $\nabla R_\eps^{\frac{1}{2}}$. Hence
\begin{equation*}
\begin{aligned}
|C_{u, \tilde u_\eta}^{(1)}| & \leq \|\bm{v}_1\|_{L^\infty} \iiint R_\eta(x - y) \Big[ (R_\eps^{\frac{1}{2}} \star \tilde u)(y - z) - (R_\eps^{\frac{1}{2}} \star \tilde u)(x - z)\Big] \ \big| \nabla R_\eps^{\frac{1}{2}}(z)\big| \tilde u(y) \ \rd z \rd y \rd x
\\
& \leq \|\bm{v}_1\|_{L^\infty} \iiint R_\eta(x) \Big[ (R_\eps^{\frac{1}{2}} \star \tilde u)(y - z) - (R_\eps^{\frac{1}{2}} \star \tilde u)(x + y - z)\Big] \ \big| \nabla R_\eps^{\frac{1}{2}}(z)\big| \tilde u(y) \ \rd z \rd y \rd x
\\
& \leq \|\bm{v}_1\|_{L^\infty} \iiint R_\eta(x) \Big[ (R_\eps^{\frac{1}{2}} \star \tilde u)(y ) - (R_\eps^{\frac{1}{2}} \star \tilde u)(x + y)\Big] \ \big| \nabla R_\eps^{\frac{1}{2}}(z)\big| \tilde u(y + z) \ \rd z \rd y \rd x,
\end{aligned}
\end{equation*}
where we use the change of variable first $x \mapsto x+y$ then $y \mapsto y + z$.
This can be further reformulated by
\begin{equation*}
\begin{aligned}
|C_{u, \tilde u_\eta}^{(1)}| \leq \|\bm{v}_1\|_{L^\infty} \iint R_\eta(x) \Big[ (R_\eps^{\frac{1}{2}} \star \tilde u)(y ) - (R_\eps^{\frac{1}{2}} \star \tilde u)(x + y)\Big] \ \big(\big| \nabla R_\eps^{\frac{1}{2}}\big| \star \tilde u\big)(y) \ \rd y \rd x,
\end{aligned}
\end{equation*}
Apply Cauchy-Schwartz inequality to the $y$ variable followed by a Sobolev estimate, we have
\begin{equation*}
\begin{aligned}
|C_{u, \tilde u_\eta}^{(1)}| & \leq \|\bm{v}_1\|_{L^\infty} \int R_\eta(x) \Big\| (R_\eps^{\frac{1}{2}} \star \tilde u) (\cdot) - (R_\eps^{\frac{1}{2}} \star \tilde u) (\cdot + x)\Big\|_{L^2} \Big\| \big| \nabla R_\eps^{\frac{1}{2}} \big| \star \tilde u\Big\|_{L^2} \ \rd x
\\
& = \|\bm{v}_1\|_{L^\infty} \bigg( \int |x| R_\eta(x) \ \rd x \bigg) \Big\| \nabla R_\eps^{\frac{1}{2}} \star \tilde u\Big\|_{L^2} \Big\| \big| \nabla R_\eps^{\frac{1}{2}} \big| \star \tilde u\Big\|_{L^2}.
\end{aligned}
\end{equation*}
Since $R$ has a bounded second moment (here, we only use the first moment), we must have that
\begin{equation*}
\begin{aligned}
\int |x| R_\eta(x) \ \rd x \leq C \eta.
\end{aligned}
\end{equation*}
Using this together with Lemma~\ref{lem:intermediate_2}, we have
\begin{equation*}
\begin{aligned}
|C_{u, \tilde u_\eta}^{(1)}| \leq C \|\bm{v}_1\|_{L^\infty} \bigg(\frac{\eta}{\eps}\bigg) \Big\| R_\eps^{\frac{1}{2}} \star \tilde u\Big\|_{L^2} \Big\| \nabla R_\eps^{\frac{1}{2}} \star \tilde u\Big\|_{L^2}.
\end{aligned}
\end{equation*}
\\

\paragraph{\emph{Bound of $C_{u, \tilde u_\eta}^{(2)}$}}
The term $C_{u, \tilde u_\eta}^{(2)}$ can be rewritten as
\begin{equation*}
\begin{aligned}
C_{u, \tilde u_\eta}^{(2)} & \defeq \int \bm{v}_1\cdot \tilde u_\eta \ (\nabla L_{\eps,\eta} \star R_\eta \star \tilde u ) \ \rd x - \int \bm{v}_1\cdot ( \tilde u_\eta \ \nabla \tilde u_\eta ) \star L_{\eps,\eta} \ \rd x
\\
& = \iint \bm{v}_1(x) \cdot L_{\eps,\eta}(x-y) (\tilde u_\eta(x) - \tilde u_\eta(y)) \ \nabla \tilde u_\eta(y) \ \rd y \rd x
\\
& = \iint \bm{v}_1(x) \cdot L_{\eps,\eta}(y) (\tilde u_\eta(x) - \tilde u_\eta(x-y)) \ \nabla \tilde u_\eta(x-y) \ \rd x \rd y
\end{aligned}
\end{equation*}
where we use the change of variable $y \mapsto x - y$.
Apply Cauchy-Schwartz in $x$ followed by a Sobolev estimate, we deduce
\begin{equation*}
\begin{aligned}
|C_{u, \tilde u_\eta}^{(2)}| & \leq \|\bm{v}_1\|_{L^\infty} \int L_{\eps,\eta}(y) \big\|\tilde u_\eta(\cdot) - \tilde u_\eta(\cdot-y)\big\|_{L^2} \big\| \nabla \tilde u_\eta \big\|_{L^2} \ \rd y
\\
& \leq \|\bm{v}_1\|_{L^\infty} \bigg( \int |y| L_{\eps,\eta}(y) \ \rd y \bigg) \big\| \nabla \tilde u_\eta \big\|_{L^2} \big\| \nabla \tilde u_\eta \big\|_{L^2}.
\end{aligned}
\end{equation*}
Recall the first moment bound property in Definition~\ref{defi:potential}, restated here
\begin{equation*}
\begin{aligned}
\int_{\R^d} |y| |L_{\eps,\eta}(y)| \ \rd y \leq C \eps.
\end{aligned}
\end{equation*}
Apply Lemma~\ref{lem:intermediate_1}, we conclude 
\begin{equation*}
\begin{aligned}
|C_{u, \tilde u_\eta}^{(2)}| & \leq C \|\bm{v}_1\|_{L^\infty} \ \eps \bigg( \frac{\eps}{\eta} \bigg)^{2k} \Big\| \nabla R_\eps^{\frac{1}{2}} \star \tilde u\Big\|_{L^2}^2.
\end{aligned}
\end{equation*}
\\

\paragraph{\emph{Bound of $C_{u, \tilde u_\eta}^{(3)}$}}
The term $C_{u, \tilde u_\eta}^{(3)}$ can be rewritten as
\begin{equation*}
\begin{aligned}
C_{u, \tilde u_\eta}^{(3)} & = \int \bm{v}_1\cdot ( \tilde u_\eta \ \nabla \tilde u_\eta ) \star L_{\eps,\eta} \ \rd x - \int ( \bm{v}_1\cdot \tilde u_\eta \ \nabla \tilde u_\eta ) \star L_{\eps,\eta} \ \rd x
\\
& = \int ( L_{\eps,\eta} * \bm{v}_1 )\cdot ( \tilde u_\eta \ \nabla \tilde u_\eta ) \ \rd x - \int \bm{v}_1\cdot \tilde u_\eta \ \nabla \tilde u_\eta \ \rd x
\\
& = \int L_{\eps,\eta} (y) (\bm{v}_1 (x - y) - \bm{v}_1 (x)) \tilde u_\eta (x) \nabla \tilde u_\eta(x) \ \rd y\rd x.
\end{aligned}
\end{equation*}
By a H\"older inequality in $x$,
\begin{equation*}
\begin{aligned}
C_{u, \tilde u_\eta}^{(3)}
& \leq \int L_{\eps,\eta} (y) \|\bm{v}_1 (\cdot - y) - \bm{v}_1 (\cdot)\|_{L^p} \| \tilde u_\eta \|_{L^q} \| \nabla \tilde u_\eta\|_{L^2} \ \rd y.
\end{aligned}
\end{equation*}
where $p,q$ are such that
\begin{equation*}
\begin{aligned}
\frac{1}{p} + \frac{1}{q} = \frac{1}{2}.
\end{aligned}
\end{equation*}
Now, take $q$ such that $H^1$ is embedded in $L^q$ (this results $p > 2$ for $d = 2$ and $p = 3$ for $d = 3$). Then, we can apply an interpolation argument to deduce
\begin{equation*}
\begin{aligned}
\|\bm{v}_1 (\cdot - y) - \bm{v}_1 (\cdot)\|_{L^p} & \leq \|\bm{v}_1 (\cdot - y) - \bm{v}_1 (\cdot)\|_{L^1}^{\frac{1}{p}} \|\bm{v}_1 (\cdot - y) - \bm{v}_1 (\cdot)\|_{L^\infty}^{1-\frac{1}{p}}
\\
& \leq 2 |y|^{\frac{1}{p}} \|\nabla \bm{v}_1\|_{TV}^{\frac{1}{p}} \|\bm{v}_1\|_{L^\infty}^{1-\frac{1}{p}}.
\end{aligned}
\end{equation*}
Apply this to our estimate of $C_{u, \tilde u_\eta}^{(3)}$, we conclude that
\begin{equation*}
\begin{aligned}
C_{u, \tilde u_\eta}^{(3)}
& \leq \bigg( \int 2 |y|^{\frac{1}{p}} L_{\eps,\eta} (y) \ \rd y \bigg) \|\nabla \bm{v}_1\|_{TV}^{\frac{1}{p}} \|\bm{v}_1\|_{L^\infty}^{1-\frac{1}{p}} \| \tilde u_\eta \|_{L^q} \| \nabla \tilde u_\eta\|_{L^2}
\\
& \leq C \|\nabla \bm{v}_1\|_{TV}^{\frac{1}{p}} \|\bm{v}_1\|_{L^\infty}^{1-\frac{1}{p}} \ \eps^{\frac{1}{p}} \|\tilde u_\eta\|_{H^1} \| \nabla \tilde u_\eta\|_{L^2}
\\
& \leq C \|\nabla \bm{v}_1\|_{TV}^{\frac{1}{p}} \|\bm{v}_1\|_{L^\infty}^{1-\frac{1}{p}} \ \eps^{\frac{1}{p}} \bigg( \frac{\eps}{\eta} \bigg)^{2k} \Big\| R_\eps^{\frac{1}{2}} \star \tilde u\Big\|_{H^1}^2,
\end{aligned}
\end{equation*}
where the last inequality is derived from Lemma~\ref{lem:intermediate_1}.
\\

As a remark, when $d = 1$, $H^1$ is embedded not only in $L^\infty$ but also in some H\"older space, hence the H\"older inequality above is not optimal.
Instead, one can take the duality
\begin{equation*}
\begin{aligned}
\|(\bm{v}_1 (\cdot - y) - \bm{v}_1 (\cdot)) \tilde u_\eta \|_{L^2} & \leq \|\bm{v}_1 (\cdot - y) - \bm{v}_1 (\cdot)\|_{W^{-1,\infty}} \| \tilde u_\eta \|_{H^1}
\\
& \leq \|\bm{v}_1 (\cdot - y) - \bm{v}_1 (\cdot)\|_{L^1} \| \tilde u_\eta \|_{H^1},
\end{aligned}
\end{equation*}
which gives us $p = 1$. In summary, the choice of $p$ is $p = d$ for $d \neq 2$ and $p > d$ for $d = 2$, corresponding to the marginal case that $H^1$ fails to be embedded in $L^\infty$.
\\

\paragraph{\emph{Conclude the proposition}}
Summarize all the bounds we have that
\begin{equation*}
\begin{aligned}
& \frac{1}{2} W_2^2(u,\tilde u_\eta) (t) - \frac{1}{2} W_2^2(u,\tilde u_\eta) (0)
\\
& \quad \leq \int_0^t - ( C_{u, \tilde u_\eta}^{(1)} + C_{u, \tilde u_\eta}^{(2)} + C_{u, \tilde u_\eta}^{(3)} ) \ \rd \tau
\\
& \quad \leq C \Bigg( \frac{\eta}{\eps} + \eps^{\frac{1}{p}} \bigg( \frac{\eps}{\eta} \bigg)^{2k} \Bigg) ( \|\nabla \bm{v}_1\|_{TV} + \|\bm{v}_1\|_{L^\infty}) \Big\| R_\eps^{\frac{1}{2}} \star \tilde u\Big\|_{L^2([0,T];H^1)}^2,
\end{aligned}
\end{equation*}
where the $\bm{v}_1$ terms are bounded by Lemma~\ref{lem:W2diff} and the $R_\eps^{\frac{1}{2}} \star \tilde u$ is bounded by the energy estimate given by Lemma~\ref{lem:energy}.
Let $\eta = \eps^{\gamma}$. Optimizing by $\gamma = 1 + \frac{1}{p(2k + 1)}$, we have
\begin{equation*}
\begin{aligned}
\frac{1}{2} W_2^2(u,\tilde u_\eta) (t) - \frac{1}{2} W_2^2(u,\tilde u_\eta) (0) \leq C \eps^{\frac{1}{p(2k + 1)}}.
\end{aligned}
\end{equation*}
Since one can choose any $p > d$, this gives the estimate in the proposition.
\end{proof}

\subsection{Conclude Theorem~\ref{thm:main}}
To conclude Theorem~\ref{thm:main} from Proposition~\ref{prop:grad_inter}, the only piece remain is to relate $W_2(u, \tilde u_\eta)$ to $W_2(u, \tilde u)$.
\begin{proof} [Proof of Theorem~\ref{thm:main}]
By triangle inequality
\begin{equation*}
\begin{aligned}
W_2(u, \tilde u)(t) & \leq W_2(u, \tilde u_\eta)(t) + W_2(\tilde u_\eta, \tilde u)(t)
\\
& \leq W_2(u, \tilde u_\eta)(0) + C \eps^{\frac{1}{p(4k + 2)}} + W_2(\tilde u_\eta, \tilde u)(t).
\end{aligned}
\end{equation*}
Note that for initial data $\tilde u_0 = u_0$. Hence
\begin{equation*}
\begin{aligned}
W_2(u, \tilde u)(t) \leq C \eps^{\frac{1}{p(4k + 2)}} + W_2(\tilde u_\eta, \tilde u)(0) + W_2(\tilde u_\eta, \tilde u)(t).
\end{aligned}
\end{equation*}
\\

Recall that $\tilde u_\eta = R_\eta \star \tilde u$. Define
\begin{equation*}
\begin{aligned}
\pi(x,y) = R_\eta(x - y) \tilde u(y).
\end{aligned}
\end{equation*}
It is easy to verify that
\begin{equation*}
\begin{aligned}
& \int \pi(x,y) \ \rd x = \tilde u(y), && \int \pi(x,y) \ \rd y = \tilde u_\eta(x)
\end{aligned}
\end{equation*}
hence $\pi$ is a transport plan between $\tilde u$ and $\tilde u_\eta$.
Moreover, one simply has
\begin{equation*}
\begin{aligned}
\int d(x,y)^2 \pi(x,y) \ \rd x \rd y = \int |x|^2 R_\eta(x) \ \rd x \ \int u(y) \ \rd y \leq C \eta^2.
\end{aligned}
\end{equation*}
Consequently
\begin{equation*}
\begin{aligned}
W_2(\tilde u_\eta, \tilde u) \leq C \eta \ll C \eps^{\frac{1}{p(4k + 2)}},
\end{aligned}
\end{equation*}
which leads us to the desired result
\begin{equation*}
\begin{aligned}
W_2(u, \tilde u) \leq C \eps^{\frac{1}{p(4k + 2)}}.
\end{aligned}
\end{equation*}

\end{proof}

\subsection{Proof of Corollary~\ref{cor:main}}
Corollary~\ref{cor:main} is derived immediately by the following:
\begin{proof}[Proof of Corollary~\ref{cor:main}]
The following interpolation argument is classical and also straightforward in the Fourier domain:
For any $k > 0$ and $f \in L^2([0,T] \times \T^d)$,
\begin{equation*}
\begin{aligned}
\big\| f \big\|_{L^2([0,T] \times \T^d)} \leq \big\| f \big\|_{L^2([0,T]; H^{-k}(\T^d))}^{\frac{1}{1+k}} \big\| f \big\|_{L^2([0,T]; H^1(\T^d))}^{\frac{k}{1+k}}.
\end{aligned}
\end{equation*}
Due to the Sobolev embedding $W^{1,-1} \subset H^{-k}$ when $k = 1 + \frac{d}{2}$, one further have
\begin{equation*}
\begin{aligned}
\big\| f \big\|_{L^2([0,T] \times \T^d)} \leq C \big\| f \big\|_{L^2([0,T]; W^{1,-1}(\T^d))}^{\frac{2}{2+d}} \big\| f \big\|_{L^2([0,T]; H^1(\T^d))}^{\frac{d}{2+d}}.
\end{aligned}
\end{equation*}
Now apply this estimate to $f = u - R_\eps^{\frac{1}{2}} \star \tilde u$. The second factor is bounded, as an a priori $L^2_tH^1_x$-energy bound independent to $\eps > 0$ has been established in Lemma~\ref{lem:energy}.
Since $f$ is now a difference of $u , R_\eps^{\frac{1}{2}} \star \tilde u \in \Pb(\T^d)$, it is well-known that the $W^{1,-1}$ norm is equivalent to the Wasserstein-$1$ distance (which are both the dual of $W^{1,\infty}$).
Therefore,
\begin{equation*}
\begin{aligned}
& \big\| u - R_\eps^{\frac{1}{2}} \star \tilde u \big\|_{L^2([0,T] \times \T^d)}
\\
& \quad \leq C \big\| u - R_\eps^{\frac{1}{2}} \star \tilde u \big\|_{L^2([0,T]; W^{1,-1}(\T^d))}^{\frac{2}{2+d}} \big\| u - R_\eps^{\frac{1}{2}} \star \tilde u \big\|_{L^2([0,T]; H^1(\T^d))}^{\frac{d}{2+d}}
\\
& \quad \leq C \bigg( \int \Big( W_1\big( u , R_\eps^{\frac{1}{2}} \star \tilde u \big) \Big)^2 \ \rd t \bigg)^{\frac{1}{2+d}}
\\
& \quad \leq C \bigg( \int 2 \Big( W_1\big( u , \tilde u \big) \Big)^2 + C \eps^2 \ \rd t \bigg)^{\frac{1}{2+d}}.
\end{aligned}
\end{equation*}
The Wasserstein-$1$ distance is further bounded by the Wasserstein-$2$ distance. Hence by Theorem~\ref{thm:main},
\begin{equation*}
\begin{aligned}
\big\| u - R_\eps^{\frac{1}{2}} \star \tilde u \big\|_{L^2([0,T] \times \T^d)} \leq C \eps^{\frac{2}{2 + d} r}.
\end{aligned}
\end{equation*}

\end{proof}

\nocite{*}
\bibliography{porous}{} 
\bibliographystyle{siam} 

\end{document}